\documentclass[a4paper,10pt]{article}
\usepackage[centertags]{amsmath}
\usepackage{amsfonts}
\usepackage{amssymb}
\usepackage{amsthm}
\usepackage{epsfig}
\usepackage{setspace}
\usepackage{ae}
\usepackage{eucal}
\usepackage[usenames]{color}%
\usepackage{graphicx}

\theoremstyle{plain}
\newtheorem{thm}{Theorem}[section]
\newtheorem{lem}[thm]{Lemma}

\newtheorem{prop}[thm]{Proposition}
\theoremstyle{definition}
\newtheorem{defi}[thm]{Definition}

\theoremstyle{remark}

\newtheorem{rem}[thm]{Remark}

\title{On global regular solution branches and multiple solutions of the Boltzmann equation}
\author{J\"org Kampen }

\begin{document}

\maketitle

\begin{abstract}
Existence of global regular solution branches of the Boltzmann Cauchy problem with continuously differentiable data in phase space dimension $2d\geq 6$ with polynomial decay at infinity of order greater than $2d$ is proved. There are data in this class of infinite relative entropy with respect to the Gaussian. Furthermore, there are weakly singular solution branches of the Boltzmann equation in spatial dimension $d\geq 3$, i.e.,  solutions of the Boltzmann equations which are only Lipschitz with respect to the velocity variables at some point in phase space. This is in accordance with a.e. $L^1$-uniqueness of renormalized solutions (cf.\cite{L}) and more classical results in function spaces of mixed regularity.
\end{abstract}

\section{Introduction}

The Boltzmann equation
\begin{equation}\label{boltz}
\partial_t F+v\nabla_xF=Q^S(F,F)
\end{equation}
determines the dynamics of the number density
\begin{equation}
{\mathbb R}\times {\mathbb R}^3\times {\mathbb R}^3\ni (t,x,v)\rightarrow F(t,x,v)\geq 0 
\end{equation}
in phase space. The right side $Q^S(F,F)$ is the Boltzmann collision integral, and, as the Boltzmann approach considers only binary collisions, it is naturally described by a bilinear operator with a characteristic dependence on the velocity $v=(v_1,v_2,v_3)$.   For a hard sphere gas the collision integral is defined by
\begin{flushleft}
\begin{equation}\label{source1}
Q^S(F,F)=\int_{(v_*,\sigma)\in {\mathbb R}^3\times S^2}
 \left( F\left(\tilde{v} \right)F\left(\tilde{v}_* \right)-
 F\left(v \right)F\left(v_* \right)\right) 
|v-v_*|dv_*d\sigma 
\end{equation}
along with 
\begin{equation}
\begin{array}{ll}
\tilde{v}=\frac{1}{2}\left(v+v_* \right) +\frac{1}{2}|v-v_*|\sigma \\
\\
\tilde{v}_*=\frac{1}{2}\left(v+v_* \right) -\frac{1}{2}|v-v_*|\sigma .
\end{array}
\end{equation}
\end{flushleft}
The results of this paper can be generalised to vector-valued mixture equations, but we consider this standard scalar equation for simplicity. Note that in the definition of the source term in (\ref{source1}) the dependence of the density function $F$ of space and time is suppressed as usual. 

We compare the following global existence theory with the global existence theory of Lion and Di Perna, which is a construction in weak function spaces.
In this global existence theory renormalized solutions are considered. These are solutions which have a certain distance to the Gaussian-Maxwell normal distribution. The required distance to this Gaussian is finite relative entropy. We recall
\begin{defi}
Given measurable functions $N,M:{\mathbb R}^3\times {\mathbb R}^3\rightarrow {\mathbb R}$ with $N\geq 0$ and $M>0$, the relative entropy of $N$ relative to $M$ is
\begin{equation}
H\left(N|M\right)=\int_{{\mathbb R}^3\times {\mathbb R}^3}\left[ N\ln\left(\frac{N}{M}\right)-N+M  \right] dxdv >0.
\end{equation}
\end{defi}
\begin{defi}
A renormalized solution of the Boltzmann equation in (\ref{boltz}) is a function
\begin{equation}
0\leq F\in C\left({\mathbb R}_+,L^1_{loc}\left({\mathbb R}^3\times {\mathbb R}^3 \right)  \right), 
\end{equation}
such that 
\begin{equation}
\frac{Q^S\left(F,F \right)}{1+F}\in  L^1_{loc}\left({\mathbb R}_+,{\mathbb R}^3\times {\mathbb R}^3 \right),
\end{equation}
and
\begin{equation}
\left(\partial_t+v\cdot\nabla_x \right)\ln\left(1+F\right)=\frac{Q^S\left(F,F \right)}{1+F} 
\end{equation}
in the sense of distributions.
\end{defi}
In the follwing we denote the parametrized Gaussian by
\begin{equation}
G_{(a,b,c)}=\frac{a}{(2\pi c)^{3/2}}\exp\left(-\frac{|x-b|^2}{2c} \right).
\end{equation}
A main result is (cf. \cite{PL})
\begin{thm}
Given initial data $F^0\equiv F(0,.)\geq 0$ a.e., where
\begin{equation}
H\left(F^0|G_{(1,0,1)} \right)<\infty,
\end{equation}
there exists a renormalized solution $F$ of the Boltzmann equation with initial data $F^0$. This renormalized solution satisfies 
\begin{equation}
\forall t\geq 0~~H\left(F(t,.)|G_{(1,0,1)} \right)+\int_0^t\int_{{\mathbb R}^3} R(F)(s,x)dxds\leq H(F(0,.)|G_{(1,0,1)})
\end{equation}
Here
\begin{equation}
R(F)=-\int_{{\mathbb R}^3}Q^S(F,F)\ln(F)dv \geq 0
\end{equation}
represents the entropy production rate.
\end{thm}
Subsequent research is usually in this framework in the sense that existence and properties of solutions with finite relative entropy with respect to the Gaussian are considered. There is an important progress concerning regularity of solution, where the connection to the Landau equation was a mayor step (cf. \cite{GS} and references therein).     
The assumption of finite relative entropy of the data with respect to the Gaussian is in some sense a strong assumption and in some sense a weak assumption. It is a weak assumption in the sense that the data may be quite irregular locally (essentially $L^1$), but it is strong in the sense that there are data which lie in strong function spaces but have no finite relative entropy distance relative to the Gaussian.
For example consider for $7<s<8$ the function $p:{\mathbb R}^3\times {\mathbb R}^3\rightarrow {\mathbb R}$ with
\begin{equation}
p(x,v):=\frac{1}{1+r^{s}}. \mbox{where}~r=\sqrt{x_1^2+x_2^2+x_3^3+v_1^2+v_2^2+v_3^2}.
\end{equation}
We have
\begin{equation}
p\in H^m\left({\mathbb R}^3\times {\mathbb R}^3 \right)~\mbox{forall positive integers $m$},
\end{equation}
but
\begin{equation}
\begin{array}{ll}
H\left(p|G_{(1,0,1)}\right)=\int_{{\mathbb R}^3\times {\mathbb R}^3}\left[ p\ln\left(\frac{p}{G_{(1,0,1)}}\right)-p+G  \right] dxdv\\
\\
=\int_{{\mathbb R}^3\times {\mathbb R}^3}\left[ 
p\left( \ln\left( p (2\pi)^{3/2}\exp\left(\frac{|x|^2}{2} \right)\right)  -1\right)+\frac{1}{(2\pi )^{3/2}}\exp\left(-\frac{|x|^2}{2} \right)  \right] dxdv\\
\\
=\int\int_{r>0}\left[ 
\frac{\left( \ln\left( \frac{1}{1+r^s} (2\pi)^{3/2}\exp\left(\frac{|r|^2}{2} \right)\right)  -1\right)}{1+r^s}+\frac{1}{(2\pi )^{3/2}}\exp\left(-\frac{|r|^2}{2} \right)  \right] r^{5}drd\Phi\\
\\
\geq\int\int_{0}^c\left[ 
\frac{ \ln\left( \frac{1}{1+r^s} \exp\left(\frac{|r|^2}{2} \right)\right)  }{1+r^s}  \right] r^{5}drd\Phi\uparrow  \infty~,\mbox{as $c \uparrow \infty$, where }~s\in (6,8),
\end{array}
\end{equation}
and where $\Phi$ is the angle function of polar coordinates.
The latter estimate indicates that the finite relative entropy condition with respect to the Gaussian is a stronger condition with respect to the required order of decay at spatial infinity than the multiplicative properties of the nonlinear part of the operator and the linear coefficients of the lineear part of the operator indicate (cf. below).  This motivates  a different perspective in this paper: we use stronger assumptions with respect to the local regularity of the data,  weaker conditions  with respect to the polynomial decay at infinity, and prove the existence of strong global solutions.
Natural function spaces for the construction of global regular solution branches of the Boltzmann equation in $2d$-dimensional phase space are
\begin{equation}\label{funct}
C^k\cap {\cal C}^{s,d}_{pol,k},~ k\geq 1,~ s> 2d.
\end{equation}
The classical phase space is $2d=6$, where$d$ denotes the spatial dimension, and we formulate our results in this case, although the generalisation to $2d\geq 6$ is straightforward.
In (\ref{funct}) $C^k$ is the function space of $k$-times continuously differentiable functions, and
\begin{equation}\label{pol1}
{\cal C}^{s,2d}_{pol,m}={\Big \{} f:{\mathbb R}^{2d}\rightarrow {\mathbb R}: \\
\\
\exists c>0~\forall |z|\geq 1~\forall 0\leq |\gamma|\leq m~{\big |}D^{\gamma}_xf(z){\big |}\leq \frac{c}{1+|z|^s} {\Big \}},
\end{equation}
where $d\geq 1$ is an arbitrary integer dimension number of space. A closed ball of 'radius' $c$ in ${\cal C}^{s,2d}_{pol,m}$ is denoted by
\begin{equation}\label{pol2}
{\cal C}^{s,2d,c}_{pol,m}={\Big \{} f:{\mathbb R}^{2d}\rightarrow {\mathbb R}: \\
\\
\forall |x|\geq 1~\forall 0\leq |\gamma|\leq m~{\big |}D^{\gamma}_xf(z){\big |}\leq \frac{c}{1+|z|^s} {\Big \}}.
\end{equation}
We construct classical solutions $F$ of the Boltzmann equation in (closed subspaces of) the function space
\begin{equation}\label{polst}
C^1\left( [0,T],C^1\cap {\cal C}^{s,2d}_{pol,1}\right) ,~\mbox{for some}~s> 2d.
\end{equation}
In the case of spatial dimension $d=3$ this means that global regular solution branches can be realised in case of polynomial decay of order $s>6$ in phase space of the initial density via compactness arguments. In this context we note that sequences in is a closed ball of the function space (such as in (\ref{pol1}) for fixed $c>0$) are directly related to sequences  classical Banach spaces on bounded domains: for $f\in {\cal C}^{2m,d}_{pol,m}$ and $x=(x_1,\cdots,x_d),~ y=(y_1,\cdots,y_d)$ the function
\begin{equation}\label{omega}
g:\Omega:=\left]-\frac{\pi}{2},\frac{\pi}{2} \right[^d\rightarrow {\mathbb R},~ g(y)=f(x),~ y_j=\arctan(x_j),~1\leq j\leq d
\end{equation} 
lives in the Banach function space $C^m\left(\left]-\frac{\pi}{2},\frac{\pi}{2} \right[^d\right)$ defined below in (\ref{cm}).  We recall
\begin{prop}
For open and bounded $\Omega\subset {\mathbb R}^n$ and consider the function space 
\begin{equation}\label{cm}
\begin{array}{ll}
C^m\left(\Omega\right):={\Big \{} f:\Omega \rightarrow {\mathbb R}|~\partial^{\alpha}f \mbox{ exists~for~}~|\alpha|\leq m\\
\\
\mbox{ and }\partial^{\alpha}f \mbox{ has an continuous extension to } \overline{\Omega}{\Big \}}
\end{array}
\end{equation}
where $\alpha=(\alpha_1,\cdots ,\alpha_n)$ denotes a multiindex and $\partial^{\alpha}$ denote partial derivatives with respect to this multiindex. Then the function space $C^m\left(\overline{\Omega}\right)$ with the norm
\begin{equation}
|f|_m:=|f|_{C^m\left(\overline{\Omega}\right) }:=\sum_{|\alpha|\leq m}{\big |}\partial^{\alpha}f{\big |}
\end{equation}
is a Banach space. Here,
\begin{equation}
{\big |}f{\big |}:=\sup_{x\in \Omega}|f(x)|.
\end{equation}

\end{prop}  
The construction is essentially in such spaces.
We prove
\begin{thm}\label{mainthm}
For initial data $F^0\equiv F(0,.)\geq 0$ a.e., where
\begin{equation}
F^0\in  C^1\cap {\cal C}^{s,6}_{pol,1},~s> 6,
\end{equation}
there exists a global classical regular solution $F$ of the Boltzmann equation with initial data $F^0$, which satisfies 
\begin{equation}
\mbox{for all $t\geq 0$:}~ F(t,.)\in  C^1\cap {\cal C}^{s-1,6}_{pol,1}
\end{equation}
Here, by 'classical solution' we mean a solution of a differential equation which satisfies the equation pointwise with classical derivatives in the sense of Weierstrass. 
\end{thm}

We note that a physical solution should satisfy $F\geq 0$ a.e., and this relation is satisfied in case of a physical classical solution, which is constructed here. If we use the term 'global regular solution branch', then we refer to our expectation that global solutions of the Boltzmann Cauchy problem are not unique in general and in case of dimension $n\geq 3$. We have
\begin{thm}\label{mainthm2}
In the situation of Theorem \ref{mainthm} and for $d\geq 3$ there exist weakly singular solution branches next to a global regular solution branch. Here the term 'weakly singular solution' means that there exists a Lipschitz continuous distributional solution  which is not classical everywhere. 
\end{thm}

\begin{rem}
 A strong form of a weakly singular solution is a solution which is classical on a domain except at one point of the phase space where it is only Lipschitz. Such solution exist. 
\end{rem}
Global regular solution branches and weakly singular solutions are constructed via viscosity limits ($\nu\downarrow 0$) of  local time solutions of the extended equations
\begin{equation}\label{boltza}
\partial_t F^{\nu}+\nu\Delta F^{\nu}+v\nabla_xF^{\nu}=Q^S(F^{\nu},F^{\nu}),~ F^{\nu}(0,.)=F^0,
\end{equation}
for the approximative particle densities $F^{\nu}(t,x,v),~ \nu>0$. A local time solution has the 
representation
\begin{equation}\label{prob2}
F^{\nu}=F^0\ast^g_{sp}\Gamma^v_{\nu}+Q^S(F^{\nu},F^{\nu})\ast^g \Gamma^v_{\nu},
\end{equation} 
where $\partial_t \Gamma^v_{\nu}+\nu\Delta \Gamma^v_{\nu}+v\nabla_x\Gamma^v_{\nu}=0$. Here, the upper script $'g'$ in $\ast^g$, and $\ast^g_{sp}$ indicates that these symbols denote generalized convoluted integrals (not usual convolutions as coefficients $v_i$  are not constant) , i.e., we define
\begin{equation}\label{astdef}
\begin{array}{ll}
F^0\ast^g_{sp}\Gamma^v_{\nu}(t,x,v):=\int_{{\mathbb R}^{2d}} F^0(y)\Gamma^v_{\nu}(t,x,v;0,y)dy\\
\\
Q^S(F^{\nu},F^{\nu})\ast^g \Gamma^v_{\nu}(t,x,v):=\int_0^t\int_{{\mathbb R}^{2d}}Q^S(F^{\nu},F^{\nu})(s,y) \Gamma^v_{\nu}(t,x,v;s,y)dyds.
\end{array}
\end{equation}
For data $F^0$ as in Theorem \ref{mainthm} or in Theorem \ref{mainthm2} we shall use properties of the adjoint $\Gamma^{v,*}_{\nu}$ of the fundamental solution $\Gamma^v_{\nu}$ in order to observe that for all $t>0$ and $z=(x,v)\in {\mathbb R}^{2d}$ there exists a $\delta \in (0,1)$ and a finite constant $C>0$ (independent of $\nu$) such that
\begin{equation}\label{visclim}
\begin{array}{ll}
{\big |}\nu\Delta F^{\nu}(t,z){\big |}={\big |}\nu\Delta \left( F^0\ast^g_{sp}\Gamma^v_{\nu}+Q^S(F^{\nu},F^{\nu})\ast^g \Gamma^v_{\nu}\right)(t,z){\big |}\\
\\
\leq {\big |}\nu\sum_i \left( F^0_{,i}\ast^g_{sp}\Gamma^{v,*}_{\nu,i}+
\left( Q^S(F^{\nu},F^{\nu})\right)_{,i} \ast^g \Gamma^{v.*}_{\nu,i}\right)(t,z){\big |}\\
\\
\leq \nu \frac{C}{(\nu t)^{\delta}}+\nu \frac{C}{(\nu)^{\delta}}\downarrow 0~ \mbox{as}~ \nu\downarrow 0~\mbox{and}~ t>0.
\end{array}
\end{equation}
Here, $C$ is a finite constant which is independent of $\nu$. The properties of the fundamental solution $\Gamma^v_{\nu}$ and its adjoint are discussed in the appendix. The limit in (\ref{visclim}) holds for data
$F^0(.)\in  C^1\cap {\cal C}^{s,6}_{pol,1}$ for $s>6$, but it holds, e.g., also for data $x\rightarrow r^{1+\alpha_0}\sin\left( \frac{1}{r^{\alpha_0}}\right)\frac{1}{1+r^8}$, where $r=\sqrt{\sum_{i=1}^dx_i^2+\sum_{i=1}^dv_i^2}$, and  $\alpha_0\in (0,5)$ which are  only Lipschitz at one point in phase space. It holds even for some examples of data which are strongly H\"older continuous at one point in space time.  
The multiplicative property the nonlinear term implies that the functional increment
\begin{equation}
\delta F^{\nu}:=F^{\nu}-F^0\ast^g_{sp}\Gamma^v_{\nu}:=Q^S(F^{\nu},F^{\nu})\ast^g \Gamma^v_{\nu}
\end{equation}
can be obtained by local time iterative solution schemes $F^{\nu,k}=F^0\ast^g_{sp}\Gamma^v_{\nu}+\delta F^{\nu,k}$ which satisfy
 \begin{equation}
\delta F^{\nu,k}:=Q^S(F^{\nu,k-1},F^{\nu,k-1})\ast^g \Gamma^v_{\nu}\in  C^1\cap {\cal C}^{6}_{pol,1},~k\geq 1,~ F^{\nu,0}=F^0\ast^g_{sp}\Gamma^v_{\nu}.
\end{equation}
It follows that the sequence $(\delta F^{\nu,*,k})_k$ with $\delta F^{\nu,*,k}(t,z)=\delta F^{\nu,k}(t,\tan(z))$ has a fixed point limit
\begin{equation}
\delta F^{\nu,*}(t,.)\in C^1\left(\Omega \right),
\end{equation}
where $\Omega$ is the bounded domain of (\ref{omega}) and $t$ is in some short time interval 
$[t_0,t_0+\Delta]$ which implies that $\delta F^{\nu}(t,.)\in C^1\cap H^1\left(\Omega \right)$.
 The function  $F^{\nu}:=F^0\ast^g_{sp}\Gamma^v_{\nu}+\delta F^{\nu}$
 satisfies the Cauchy problem for all $\nu >0$ and the corresponding family $F^{\nu,*}:=(F^{0}\ast^g_{sp}\Gamma^v_{\nu})(.,\tan(.))+\delta F^{\nu}(.,\tan(.))$ turns out to be a pointwise bounded and equicontinuous family of functions with a bounded subsequence $F^{\nu_k,*}, \nu_k\downarrow 0$ which converges along with $F^{\nu_k}, \nu_k\downarrow 0$ to $F^*$ and $F$ respectively with $F\in C\left(\left[t_0,t_0+\Delta \right],C^1\cap H^1\right)$. Moreover, the increment $\delta F^{\nu}$ inherits the order of polynomial spatial decay of order $s>6$, where here and henceforth by 'spatial polynomial decay' we mean decay with respect to the phase space variables as their modulus goes to infinity .  Since (\ref{visclim}) holds the function $F$ is a classical local time solution of the Boltzmann equation. Using the structure of the nonlinear term we shall observe in the next section that this local time solution can be extended to a global classical solution branch. As the local time argument can be performed for some mentioned data which are only Lipschitz at one point in phase space argument it follows that Boltzmann's gas theory is not deterministic globally, or intrinsically incomplete. This situation has some similarities with the case of the incompressible Euler equation, since both equations do not loose  well-posedness under the transformation $t\rightarrow -t$ (although we have a characterisation of a time direction by the entropy in case of the Boltzmann equation). This is in accordance with Lions result that there is a.e. uniqueness of weak normalized solutions. In the next we provide the details of this proof. In the appendix we add some information concerning the fundamental solution of $\Gamma^v_{\nu}$, its adjoint $\Gamma^{v,*}_{\nu}$, and a priori estimates of these distributions  and their spatial derivatives. Most of this material is known from stochastic analysis.
\section{Proof of Theorem \ref{mainthm}}

For the extended Cauchy problem in (\ref{boltza}) in item i) below we prove local time existence of solutions $F^{\nu}$ in $C^1\left( \left[t_0,t_0+\Delta \right],H^1\cap C^1\right) $ for $t_0\geq 0$ and small $\Delta >0$ and for data  in $F^{\nu,t_0}=F(t_0,.)\in  C^1\cap {\cal C}^{s,6}_{pol,1}$. Then in item ii) below we obtain a viscosity limit $F^{\nu_k}\rightarrow F$ for some zero sequence $(\nu_k)_k$, i.e. $\lim_{k\uparrow \infty} \nu_k=0$, which is a short time solution of the original Boltzmann equation. Furthermore, we show that the increment $\delta F=\lim_{k\uparrow \infty}\delta F^{\nu_k}$ and first order derivatives with respect to phase space variables preserve spatial polynomial decay of order $s>6$.  Finally, in item iii) we use the special structure of the nonlinear term in order to show that the increment of the nonlinear term is small compared to potential damping induced by a simple time transformation such that a global regular upper bound can be constructed which is linear in time. 

\begin{itemize}
\item[i)] Local solutions of the Boltzmann Cauchy problem on a time interval $[t_0,t_0+\Delta]$ with data $F^{t_0}\in  C^1\cap {\cal C}^{s,6}_{pol,1}$ are  constructed by a viscosity limit (positive real parameter $\nu \downarrow 0$) of the extended system
\begin{equation}\label{boltz}
\partial_t F^{\nu}-\nu\Delta F^{\nu}+v\nabla_xF^{\nu}=Q^S(F^{\nu},F^{\nu}),~ F^{\nu}(0,.)=F^0,
\end{equation}
to be solved for the approximative particle density $(t,x,v)\rightarrow F^{\nu}(t,x,v)$. Here $\partial_t F^{\nu}=F^{\nu}_{,t}$ denotes the partial derivative with respect to the time variable, and $x,v\in {\mathbb R}^d$. If the coefficients $v_i,~ 1\leq i\leq d$ have a linear upper bound (which will be verified by a local iteration scheme), then the fundamental solution
\begin{equation}\label{gausext}
\partial_t \Gamma^v_{\nu}+\nu\Delta \Gamma^v_{\nu}+v\nabla_x\Gamma^v_{\nu}=0.
\end{equation}
exists, and a  local solution of the extended Boltzmann Cauchy problem on a time interval $[t_0,t_0+\Delta]$ has the representation (note that $\ast^g,\ast^g_{sp}$ are {\it not} standard convolutions, cf.  the definition in (\ref{astdef}) above)
\begin{equation}\label{prob2}
F^{\nu}=F^0\ast^g_{sp}\Gamma^v_{\nu}+Q^S(F^{\nu},F^{\nu})\ast^g \Gamma^v_{\nu},~~~\mbox{(cf.~\ref{astdef})}.
\end{equation} 
An alternative representation of the local solution is of the form
\begin{equation}\label{prob1}
F^{\nu}=F^0\ast_{sp}G_{\nu}-v\nabla_xF^{\nu}\ast G_{\nu}+Q^S(F^{\nu},F^{\nu})\ast G_{\nu},
\end{equation}
where $G_{\nu}$ is the fundamental solution of the equation $\partial_tG_{\nu}-\nu \Delta G_{\nu}=0$. Note that in the latter case $\ast$ and $\ast_{sp}$ denote indeed standard convolutions.
For positive $\nu >0$ we solve these fixed point equations by an iteration scheme and then consider the viscosity limit.
A local iteration scheme for the Boltzmann Cauchy problem on the time interval $[t_0,t_0+\Delta]$ is given by the list of functions $(F^{\nu}_k)_{k\geq 0}$, where $F^{\nu}_0:=F^0\ast^g_{sp}\Gamma^v_{\nu}$ and for $k \geq 1$ the approximation $F^{\nu}_k$ is determined recursively by
\begin{equation}\label{prob1iter}
F^{\nu}_k=F^0\ast^g_{sp}\Gamma^v_{\nu}+Q^S(F^{\nu}_{k-1},F^{\nu}_{k-1})\ast^g \Gamma^v_{\nu}.
\end{equation}
For the first functional increment of this series we have 
\begin{equation}
\begin{array}{ll}
\delta F^{\nu}_1:=F^{\nu}_1-F^0\ast^g_{sp}\Gamma^v_{\nu}
=Q^S(F^{\nu}_{0},F^{\nu}_{0})\ast^g \Gamma^v_{\nu},
\end{array}
\end{equation}
and for $k\geq 1$ we have
\begin{equation}
\begin{array}{ll}
\delta F^{\nu}_{k+1}:=F^{\nu}_{k+1}-F^{\nu}_{k}=\delta Q^S(F^{\nu}_{k},F^{\nu}_{k})\ast^g \Gamma^v_{\nu},
\end{array}
\end{equation}
where
\begin{equation}
\begin{array}{ll}
\delta Q^S(F^{\nu}_{k},F^{\nu}_{k})\ast^g \Gamma^v_{\nu}=\left( Q^S(F^{\nu}_{k},F^{\nu}_{k})-Q^S(F^{\nu}_{k-1},F^{\nu}_{k-1})\right) \ast^g \Gamma^v_{\nu}.
\end{array}
\end{equation}
We write 
\begin{equation}\label{deltaQ}
\begin{array}{ll}
\delta Q^S(F^{\nu}_{k},F^{\nu}_{k})= Q^S(F^{\nu}_{k},F^{\nu}_{k})-Q^S(F^{\nu}_{k-1},F^{\nu}_{k})\\
\\
\hspace{2.2cm}+Q^S(F^{\nu}_{k-1},F^{\nu}_{k})-Q^S(F^{\nu}_{k-1},F^{\nu}_{k-1})
\end{array}
\end{equation}
For functions $H:[t_0,t_{0}+\Delta]\times {\mathbb R}^d\times {\mathbb R}^d\rightarrow {\mathbb R}$ we define
\begin{equation}\label{nH00}
\begin{array}{ll}
{\big |}H{\big |}^{t_0,\Delta}_{sup}:=
\sup_{(t,x,v)\in [t_0,t_0+\Delta] \times {\mathbb R}^d\times {\mathbb R}^d}
 {\big |}H(t,x,v)  {\big |},
\end{array}
\end{equation}
and
\begin{equation}\label{nH0}
\begin{array}{ll}
 {\big |}H{\big |}^{t_0,\Delta}_{sup,1}:=
\sup_{(t,x,v)\in [t_0,t_0+\Delta] \times {\mathbb R}^d\times {\mathbb R}^d}\left( {\big |}H(t,x,v) {\big |}+{\big |}\nabla_{(x,v)}H(t,x,v) {\big |}\right),
\end{array}
\end{equation}
where $\nabla_{(x,v)}$ denotes the gradient with respect to the spatial phase spaces variables$(x,v)\in {\mathbb R}^{2d}$. We also write
\begin{equation}\label{nH1}
{\big |}H{\big |}^{t_0,\Delta,s}_{sup}:={\big |}(1+r^s)H{\big |}^{t_0,\Delta}_{sup},
\end{equation}
and
\begin{equation}\label{nH2}
{\big |}H{\big |}^{t_0,\Delta,s}_{sup,1}:={\big |}(1+r^s)H{\big |}^{t_0,\Delta}_{sup,1},
\end{equation}
where $r=\sqrt{\sum_{i=1}^dx_i^2+\sum_{i=1}^dv_i^2}$ is the radial variable of the phase space.
Assume inductively that for some finite constant $C>0$ we have
\begin{equation}
\begin{array}{ll}
 {\big |}F^{\nu}_{k-1}{\big |}^{t_0,\Delta,s-1}_{sup,1}\leq C,
\end{array}
\end{equation}
where we account for the linear growth of the first order term.
We continue to suppress the dependence on $(t,x)$ for brevity. 
For the first difference on the right side of (\ref{deltaQ}) we have
\begin{equation}\label{source12}
\begin{array}{ll}
{\big |}Q^S(F^{\nu}_{k},F^{\nu}_{k})(v)-Q^S(F^{\nu}_{k-1},F^{\nu}_{k})(v){\big |}=\\
\\
{\big |}\int_{(v_*,\sigma)\in {\mathbb R}^3\times S^2}
 \left(\delta F^{\nu}_k\left(\tilde{v} \right)F^{\nu}_{k}\left(\tilde{v}_* \right)-
 \delta F^{\nu}_k\left(v \right) F^{\nu}_k\left(v_* \right)\right) 
|v-v_*|dv_*d\sigma{\big |}\\
\\
\leq {\Big |}\int_{r>0}
 \left(\frac{c{\big |}\delta F^{\nu}_k{\big |}^{t_0,\Delta}_{sup}}{1+r^{s}}\frac{C}{1+r^{s-1}}\right)rr^{d-1}dr{\Big |}\leq \frac{\tilde{C}{\big |}\delta F^{\nu}_k{\big |}^{t_0,\Delta}_{sup}}{1+r^{2s-d-2}}
\end{array}
\end{equation}
for some finite constant $c$ which depends only on phase space dimension $2d$.
Similarly, for the second difference on the right side of (\ref{deltaQ}) we have
\begin{equation}\label{source12}
\begin{array}{ll}
{\big |}Q^S(F^{\nu}_{k},F^{\nu}_{k})(v)-Q^S(F^{\nu}_{k-1},F^{\nu}_{k})(v){\big |}=\\
\\
{\big |}\int_{(v_*,\sigma)\in {\mathbb R}^3\times S^2}
 \left( F^{\nu}_k\left(\tilde{v} \right)\delta F^{\nu}_{k}\left(\tilde{v}_* \right)-
 \ F^{\nu}_k\left(v \right) \delta F^{\nu}_k\left(v_* \right)\right) 
|v-v_*|dv_*d\sigma{\big |}\\
\\
\leq {\Big |}\int_{r>0}
 \left(\frac{c{\big |}\delta F^{\nu}_k{\big |}^{t_0,\Delta}_{sup}}{1+r^{s}}\frac{C}{1+r^{s-1}}\right)rr^{d-1}dr{\Big |}\leq \frac{\tilde{C}{\big |}\delta F^{\nu}_k{\big |}^{t_0,\Delta}_{sup}}{1+r^{2s-d-2}}.
\end{array}
\end{equation}
Hence, using a priori estimates for $\Gamma^v_{\nu}$ and first order phase space derivatives of $\Gamma^v_{\nu}$  for $\Delta >0$ small enough we have
\begin{equation}
\begin{array}{ll}
{\big |}\delta F^{\nu}_{k+1}{\big |}^{t_0,\Delta,s}_{sup,1}={\big |}\delta Q^S(F^{\nu}_{k},F^{\nu}_{k})\ast^g \Gamma^v_{\nu}{\big |}^{t_0,\Delta,s}_{sup,1}\leq\frac{1}{2}{\big |}\delta F^{\nu}_{k}{\big |}^{t_0,\Delta,s}_{sup,1},
\end{array}
\end{equation}
such that we have the desired contraction. 
We remark that for first order ($|\alpha|=1)$ derivatives $D^{\alpha}_z$ with respect to phase space variables $z=(x,v)$ we may use the adjoint in order to obtain 
\begin{equation}
\begin{array}{ll}
D^{\alpha}_z\delta F^{\nu}_{k+1}=\left( D^{\alpha}_z\delta Q^S(F^{\nu}_{k},F^{\nu}_{k})\right) \ast^g \Gamma^{v,*}_{\nu},
\end{array}
\end{equation}
where $D^{\alpha}_z\delta Q^S(F^{\nu}_{k},F^{\nu}_{k})$ can be represented in terms of $D^{\beta}_z\delta F^{\nu}_{k}$ for $0\leq |\beta|\leq 1$. Hence, the estimates are naturally indpendent of $\nu$ as is the contraction. For the solution one order of decay at phase space infinity is lost due to the fact that the linear first order coefficient implies a linear factor for the Gaussian upper bound of the fundamental solution $\Gamma^v_{\nu}$ (cf. appendix). 
%

\item[ii)] Next we consider viscosity limits for the local time solution of item i), i.e. on a time interval where we have contraction.  For each $\nu >0$ we have obtained fixed point solutions $F^{\nu}$ of the extended Boltzmann equation, where for all multiindices $\alpha=(\alpha_1,\cdots,\alpha_{2d})$ wit $0\leq |\alpha|\leq 1$ and $z=(x,v)$ we have
\begin{equation}\label{prob2a}
D^{\alpha}_zF^{\nu}=F^0\ast^g_{sp}D^{\alpha}_z\Gamma^v_{\nu}+Q^S(F^{\nu},F^{\nu})\ast^g D^{\alpha}_z \Gamma^{v,*}_{\nu},~ t\in [t_0,t_0+\Delta],
\end{equation} 
or
\begin{equation}\label{prob1a}
D^{\alpha}_zF^{\nu}=F^0\ast_{sp}D^{\alpha}_zG_{\nu}-v\nabla_xF^{\nu}\ast D^{\alpha}_zG_{\nu}+Q^S(F^{\nu},F^{\nu})\ast D^{\alpha_z}G_{\nu},
\end{equation}
also for $t\in [t_0,t_0+\Delta]$. Since the contraction argument in item i) holds independently of the viscosity $\nu >0$ we have 
\begin{equation}
{\Big |}F^{\nu}(t,.){\Big |}^{t_0,\Delta,s-1}_{sup,1}\leq C,~ 
\end{equation}
where $C>0$ is a finite constant which is independent of $\nu >0$. Note the upper script $s-1$ due to the linear coefficient of the linear term of the Boltzmann equation. Moreover and accordingly, the fundamental solutions and their spatial derivatives  up to first order  involved in (\ref{prob2a}) have  the a prior upper bound
\begin{equation}\label{gammaalpha}
{\Big |}D^{\alpha}_z\Gamma^{v,*}_{\nu}(t,z,s,y){\Big|}\leq c\frac{(1+|v|)}{\nu^{\frac{2d}{2}}(t-s)^{\frac{2d+|\alpha|}{2}}}\exp\left(-\lambda\frac{|z-y|^2}{\nu(t-s)} \right) 
\end{equation}
for some finite constants $\lambda, c>0$, and where $0\leq |\alpha|\leq 1$. A similar upper bound without factor $1+|v|$ holds for the Gaussian $G_{\nu}$, of course. The estimates in (\ref{gammaalpha}) imply that we have local estimates
\begin{equation}\label{gammaalpha2}
{\Big |}D^{\alpha}_z\Gamma^{v,*}_{\nu}(t,z,s,y){\Big|}\leq 
\frac{\tilde{C}(1+|v|)}{\nu^{\delta}(t-s)^{\delta}|z-y|^{2d+|\alpha|-2\delta}},~z=(x,v),
\end{equation}
and where $\tilde{C}$ is a finite constant.
For $0\leq |\alpha|\leq 1$ and $z=(x,v)$ the families $\left( D^{\alpha}_z F^{\nu}\right)_{\nu >0}$ are  equicontinuous and pointwise bounded (by the constant $C$). Since $F^\nu(.)\in  C^1\cap {\cal C}^{s,6}_{pol,1}$ for some $s>6$  the corresponding families $\left( D^{\alpha}_z F^{\nu}(.,\tan(.))\right)_{\nu >0}$ for $0\leq |\alpha|\leq 1$ and $z=(x,v)$
are uniformly pointwise bounded and equicontinuous families in $C^{1}\left(\left[t_0,t_0+\Delta\right],\Omega)  \right)$ and we get a zero sequence $(\nu_k)_k$ and limits
\begin{equation}\label{blimit1}
D^{\alpha}_z F(.,\tan(.)):=\lim_{k\uparrow \infty}D^{\alpha}_z F^{\nu_k}(.,\tan(.)),~0\leq |\alpha|\leq 1.
\end{equation} 
where we have a corresponding limit in the whole space
\begin{equation}\label{blimit2}
D^{\alpha}_z F:=\lim_{k\uparrow \infty}D^{\alpha}_z F^{\nu_k},~0\leq |\alpha|\leq 1.
\end{equation} 
Note that for the latter we have $D^{\alpha}_z F(t,.)\in {\cal C}^{s,6}_{pol,1}$ for all $t\in [t_0,t_0+\Delta]$. We remark that continuous spaces of functions on Euclidean space of any dimension which vanish at spatial infinity are indeed closed, but this may be less well-known.
Here we note also that continuous differentiability with respect to time follows from continuous differentiability with respect to variables of phase space.
In order to show that the functions in (\ref{blimit1}) and (\ref{blimit2}) are solutions of the Boltzmann Cauchy problem (in transformed coordinates in the former case of course), we first remark that all the members $F^{\nu_k},~ k\uparrow \infty$ of the series above are classical local time solutions of the extended equations
\begin{equation}\label{boltza}
\partial_t F^{\nu_k}+\nu_k\Delta F^{\nu_k}+v\nabla_xF^{\nu_k}=Q^S(F^{\nu_k},F^{\nu_k}),~ F^{\nu_k}(0,.)=F^0,~ \nu_k >0.
\end{equation}
We have to show that
\begin{equation}\label{visclim0}
\begin{array}{ll}
\lim_{\nu_k\downarrow 0}{\big |}\nu_k\Delta F^{\nu_k}(t,z){\big |}=0
\end{array}
\end{equation}
\begin{equation}\label{gammaalpha}
{\Big |}D^{\alpha}_z\Gamma^{v,*}_{\nu}(t,z,s,y){\Big|}\leq c\frac{(1+|v|)}{\nu^{\frac{2d}{2}}(t-s)^{\frac{2d+|\alpha|}{2}}}\exp\left(-\lambda\frac{|z-y|^2}{\nu(t-s)} \right) 
\end{equation}
for some finite constants $\lambda, c>0$, and where $0\leq |\alpha|\leq 1$. A similar upper bound without factor $1+|v|$ holds for the Gaussian $G_{\nu}$, of course. The estimates in (\ref{gammaalpha}) imply that we have local estimates
\begin{equation}\label{gammaalpha2}
{\Big |}D^{\alpha}_z\Gamma^{v,*}_{\nu}(t,z,s,y){\Big|}\leq 
\frac{\tilde{C}(1+|v|)}{\nu^{\delta}(t-s)^{\delta}|z-y|^{2d+|\alpha|-2\delta}},~z=(x,v),
\end{equation}

We use the adjoint of the fundamental solution and a priori estimates in (\ref{gammaalpha}) such that for some finite constants $C,C_0>0$ and for $t>0$ we have
\begin{equation}\label{visclim}
\begin{array}{ll}
{\big |}\nu\Delta F^{\nu}(t,z){\big |}={\big |}\nu\Delta \left( F^0\ast^g_{sp}\Gamma^v_{\nu}+Q^S(F^{\nu},F^{\nu})\ast^g \Gamma^v_{\nu}\right)(t,z){\big |}\\
\\
\leq {\big |}\nu\sum_{i=1}^{2d} \left( F^0_{,i}\ast^g_{sp}\Gamma^{v,*}_{\nu,i}+
\left( Q^S(F^{\nu},F^{\nu})\right)_{,i} \ast^g \Gamma^{v.*}_{\nu,i}\right)(t,z){\big |}\\
\\
\leq \int_{{\mathbb R}^{2d}}2d{\Big |}\nu \frac{C}{1+|z|^{2d}}\frac{(1+|v|)}{\nu^{\frac{2d}{2}}t)^{\frac{2d+|\alpha|}{2}}}\exp\left(-\lambda\frac{|z-y|^2}{\nu t} \right){\Big |}dz+\\
\\
\int_{t_0}^{t_0+\Delta}\int_{{\mathbb R}^{2d}}2d{\Big |}\nu \frac{C}{1+|z|^{2d}}\frac{(1+|v|)}{\nu^{\frac{2d}{2}}(t-s)^{\frac{2d+|\alpha|}{2}}}\exp\left(-\lambda\frac{|z-y|^2}{\nu(t-s)} \right){\Big |}dzds \\
\\
\leq \int_{B^{2d}_1(z)}2d{\Big |}\nu \frac{C}{1+|z|^{2d}}\frac{\tilde{C}(1+|v|)}{\nu^{\delta}(t-s)^{\delta}|z-y|^{2d+|\alpha|-2\delta}}{\Big |}dz+\nu C\\
\\
+\int_{t_0}^{t_0+\Delta}\int_{B^{2d}_1(z)}2d{\Big |}\nu \frac{C}{1+|z|^{2d}}\frac{\tilde{C}(1+|v|)}{\nu^{\delta}(t-s)^{\delta}|z-y|^{2d+|\alpha|-2\delta}}{\Big |}dzds\\
\\
\leq \nu \frac{C_0}{(\nu t)^{\delta}}+\nu \frac{C_0}{(\nu)^{\delta}}\downarrow 0~ \mbox{as}~ \nu\downarrow 0.
\end{array}
\end{equation}
Hence the viscosity limit $\nu_k\downarrow 0$ is indeed a local time solution of the Boltzmann equation.
\item[iii)]
The difference structure of the nonlinear source term $Q^S(F,F)$ encodes some spatial effects of the operator and can be exploited in order to obtain an upper bound of the solution increment $\delta F$ over a local time interval $\left[t_0.t_0+\Delta_0\right]$ which is small relative to a potential damping term introduced by a global time transformation. We start with a  classical local time representation of the value function increment $\delta F^{\nu}$  in terms of the Gaussian and first order derivatives of the Gaussian $G_{\nu}$ or of the fundamental solution  $\Gamma^{v}_{\nu}$ or its adjoint $\Gamma^{v}_{\nu}$ and their first order spatial derivatives (there are several alternatives here). First note that for any $\nu >0$  the first order spatial derivatives of the increment $\delta F^{\nu}$ we have on a time interval $t\in [t_0,t_0+\Delta]$ the natural representations  
\begin{equation}\label{prob2aiii}
D^{\alpha}_z\delta F^{\nu}=Q^S(F^{\nu},F^{\nu})\ast^g D^{\alpha}_z \Gamma^{v,*}_{\nu},~ |\alpha|=1,
\end{equation} 
or
\begin{equation}\label{prob1aiii}
D^{\alpha}_z\delta F^{\nu}=-v\nabla_xF^{\nu}\ast D^{\alpha}_zG_{\nu}+Q^S(F^{\nu},F^{\nu})\ast D^{\alpha}_zG_{\nu},~ |\alpha|=1.
\end{equation}
Concerning the increment $\delta F^{\nu}$ of the value function  we shall use the local solution representation in (\ref{prob1}) or the representation in (\ref{prob2}) and similar representations which can be derived from them. In case of convoluted representations in terms of representations with the Gaussian there is an additional flux term (compared to (\ref{prob2})). Concerning this linear drift term we remark that we can rewrite (\ref{prob1}) with
\begin{equation}
\begin{array}{ll}
v\nabla_xF^{\nu}\ast G_{\nu}=\sum_{i=1}^dv_iF^{\nu}_{,i}\ast G_{\nu}=\sum_{i=1}^dv_iF^{\nu}\ast G_{\nu,i},
\end{array}
\end{equation}  
using the convolution rule, and because the coefficients $v_i$ are not affected by the spatial derivatives $F^{\nu}_{,i}=\partial_{x_i}F^{\nu}$. We then have Lipschitz continuity of the terms $v_iF^{\nu}$ on local time intervals which we shall use together with the symmetry of the first spatial derivatives of the Gaussian $G_{\nu,i}$ as outlined above. The nonlinear term of the Boltzmann equation is itself a difference which can be represented by a multivariate Taylor remainder term. 
More precisely, we shall determine functions $P^S_j,~ 1\leq j\leq d$ and $W^S_j,~ 1\leq j\leq 2d$ such that a  local solution has on a time interval $[t_0,t_0+\Delta]$  the representation
\begin{equation}\label{prob2aaa}
\begin{array}{ll}
\delta F^{\nu}=F^{\nu}-F^0\ast_{sp}G_{\nu}=\sum_{i=1}^dv_iF^{\nu}\ast G_{\nu,i}
+\sum_{j=1}^{2d}W^S_j(F^{\nu},F^{\nu})\ast G_{\nu}\\
\\
=\sum_{j=1}^dP^S_j(F^{\nu},F^{\nu})\ast G_{\nu,v_j}+\sum_{j=d+1}^{2d}W^S_j(F^{\nu},F^{\nu})\ast G_{\nu},
\end{array}
\end{equation} 
or, alternatively, the representation
\begin{equation}\label{prob2aaa}
\begin{array}{ll}
\delta F^{\nu}=F^{\nu}-F^0\ast^g_{sp}\Gamma^v_{\nu}\\
\\
+\sum_{j=1}^{2d}W^S_j(F^{\nu},F^{\nu})\ast^g_{sp}\Gamma^v_{\nu}\\
\\
=\sum_{j=1}^dP^S_j(F^{\nu},F^{\nu})\ast^g \Gamma^{v,*}_{\nu,v_j}+\sum_{j=d+1}^{2d}W^S_j(F^{\nu},F^{\nu})\ast^g \Gamma^v_{\nu}
\end{array}
\end{equation} 
where the symbols $\Gamma^{v,*}_{\nu,v_j}$ refer to derivatives of $\Gamma^{v,*}_{\nu}$ with respect to  the components of the velocity components $v_i$. The functions $P^S_j,~ 1\leq j\leq d$ and $W^S_j,~ 1\leq j\leq 2d$ are determined below by application a multivariate Taylor formula.
We mention that for data with more regularity higher regularity of the existence scheme may be obtained using in addition representations of higher multivariate derivatives $D^{\alpha}_{z}F^{\nu}$ with respect to the phase space variables $z=(z_1,\cdots,z_{2d})=(x,v)$, where we have for $1\leq |\alpha|=|\gamma|+1$ with $\alpha_j=\beta_j+1$ and $\beta_k=\alpha_k$ for $k\neq j$ the simple representation
\begin{equation}\label{prob2aab}
D^{\alpha}_z\delta F^{\nu}=
D^{\alpha}_zF^{\nu}-D^{\beta}_zF^0\ast^g_{sp}
\Gamma^{v}_{\nu}
=D^{\beta}_zQ^S(F^{\nu},F^{\nu})\ast^g \Gamma^{v,*}_{\nu,j}.
\end{equation} 
However, the construction of global classical solutions is essential, where for $|\alpha|=1$ and $|\beta|=0$.
We shall show that due the special structure of the nonlinear term and/or (generalised) convolutions with first order derivatives of the fundamental solution $\Gamma^{v,*}_{\nu}$ the increments in (\ref{prob2aaa}) and (\ref{prob2aab}) have upper bounds which are small enough such the growth can be offset by damping potentials introduced by auto-control transformation which are pure time transformations.
%
In order to rewrite the source term
we first recall the multivariate Taylor formula. We have
\begin{lem}
For $1\leq m\geq l\geq 0$ and $g\in C^m({\mathbb R}^n)$ we have for all $x,h\in {\mathbb R}^n$
\begin{equation}
\begin{array}{ll}
g(x+h)=g(x)+\sum_{0<\alpha <l}\frac{(\partial^{\alpha}g)(x)}{\alpha !}\\
\\
+l\sum_{|\alpha|=m}\frac{h^{\alpha}}{\alpha !}\int_0^1(1-\theta)^{l-1}\left(\partial^{\alpha}g \right)(x+\theta h)d\theta ,
\end{array}
\end{equation}
where $\partial^{\alpha}g$ denotes the multivariate partial derivative of $g$ of order $\alpha$ with respect to the multiindex $\alpha=(\alpha_1,\cdots,\alpha_n)$. Note furthermore that $h^{\alpha}=h_1^{\alpha_1}\cdot h^{\alpha_2}_2\cdots h^{\alpha_n}_n$ as usual.
\end{lem}
Next define
\begin{equation}
\begin{array}{ll}
v-v^*=:w,~dv^*=-dw.
\end{array}
\end{equation}
Then $v^*=v-w$ and
\begin{equation}
\begin{array}{ll}
\tilde{v}=\frac{1}{2}\left(v+v-w \right) +\frac{1}{2}|w|\sigma =v-\frac{w}{2}+|w|\frac{\sigma}{2}\\
\\
\tilde{v}_*=\frac{1}{2}\left(v+v-w \right) -\frac{1}{2}|w|\sigma =v-\frac{w}{2}-|w|\frac{\sigma}{2}\\
\\
\hspace{0.45cm}=v-w+\frac{w}{2}-|w|\frac{\sigma}{2}.
\end{array}
\end{equation}
Next define $G$ with
\begin{equation}
G(v,w):=F(v)F_{-}(w-v),~\mbox{where for all $v\in {\mathbb R}^3$}~F_{-}(v):=F(-v)
\end{equation}
(recall that we keep the dependence on $x\in {\mathbb R}^d$).We shall use the equality
\begin{equation}
G(v+h,w+2h)=F(v+h)F_{-}(w+h-v)
\end{equation}
for 
\begin{equation}\label{hlabel}
h=-\frac{w}{2}+\frac{\sigma}{2}|w|.
\end{equation}
The source term $Q^S(F,F)$ can be written in the form 
\begin{equation}\label{source1trans}
\begin{array}{ll}
\int_{(v_*,\sigma)\in {\mathbb R}^3\times S^2}
 \left( F\left(\tilde{v} \right)F\left(\tilde{v}_* \right)-
 F\left(v \right)F\left(v_* \right)\right) 
|v-v_*|dv_*d\sigma \\
\\
=\int_{(w,\sigma)\in {\mathbb R}^3\times S^2}
 {\Big (} F\left(v-\frac{w}{2}+\frac{\sigma}{2}|w|\right)F\left(v-w+\frac{w}{2}-\frac{\sigma}{2}|w|\right)-\\
 \\
 \hspace{3cm}F\left(v \right)F\left(v-w \right){\Big )}|w|dwd\sigma\\
 \\
=\int_{(w,\sigma)\in {\mathbb R}^3\times S^2}
 {\Big (} F\left(v-\frac{w}{2}+\frac{\sigma}{2}|w|\right)F_-\left(w-v-\frac{w}{2}+\frac{\sigma}{2}|w|\right)-\\
 \\
 \hspace{3cm}F\left(v \right)F_-\left(w-v \right){\Big )}|w|dwd\sigma\\
 \\
 =\int_{(w,\sigma)\in {\mathbb R}^3\times S^2}
 {\Big (} G\left(v+h,w+2h\right)-G\left(v,w \right){\Big )}|w|dwd\sigma .
\end{array}
\end{equation}
Assuming $f\in C^1$ we can expand $G$ by the multivariate Taylor formula to lowest order. We get with $h$ as in (\ref{hlabel})
\begin{equation}
\begin{array}{ll}
G\left(v+h,w+2h\right)-G\left(v,w \right)=\\
\\
+\sum_{i=1}^{2d}\int_0^1(1-\theta)\partial_{z'_i}G\left( v+\theta h,w+\theta 2h\right)d\theta h_i,
\end{array}
\end{equation}
where $\partial_{z'_i}=\partial_{v_i}$for $1\leq i\leq d$ and $\partial_{z'_i}=\partial_{w_i}$for $d+1\leq i\leq 2d$.   
Hence, with $h=(h_1,\cdots,h_d,h_{d+1},\cdots h_{2d})$ as in (\ref{hlabel}) we have
\begin{equation}\label{source1trans}
\begin{array}{ll}
Q^S(F,F)=\\
\\
\int_{(w,\sigma)\in {\mathbb R}^3\times S^2}
 {\Big (}\sum_{i=1}^{2d}\int_0^1(1-\theta)\partial_{z'_i}G\left( v+\theta h,w+\theta 2h\right)d\theta  h_i{\Big )}|w|dwd\sigma .
\end{array}
\end{equation}
Now we can determine the functions $W^S_j,~1\leq j\leq 2d$ and  $P^S_j,~1\leq j\leq d$ 
where we may use usual convolutions as in
\begin{equation}\label{prob2aaatrans}
\begin{array}{ll}
\delta F^{\nu}=F^{\nu}-F^0\ast G_{\nu}=v\nabla_xF^{\nu}\ast G_{\nu}+Q^S(F^{\nu},F^{\nu})\ast G_{\nu}\\
\\
\hspace{0.8cm}=v\nabla_xF^{\nu}\ast G_{\nu}+\sum_{j=1}^{2d}W_j(F^{\nu},F^{\nu})\ast G_{\nu}\\
\\
\hspace{0.8cm}=v\nabla_xF^{\nu}\ast G_{\nu}+\sum_{j=1}^dP^S_j(F^{\nu},F^{\nu})\ast G_{\nu,v_j}\\
\\
\hspace{0.8cm}+\sum_{j=1}^dW^S_j(F^{\nu},F^{\nu})\ast G_{\nu}
\end{array}
\end{equation} 
or generalised convoluted integrals and the adjoint $\Gamma^{v,*}_{\nu}$ as in
\begin{equation}\label{prob2aaatrans2}
\begin{array}{ll}
\delta F^{\nu}=F^{\nu}-F^0\ast^g_{sp}\Gamma^v_{\nu}=Q^S(F^{\nu},F^{\nu})\ast^g \Gamma^v_{\nu}\\
\\
=\sum_{j=1}^{2d}W^S_j(F^{\nu},F^{\nu})\ast^g \Gamma^v_{\nu}\\
\\
\\
=\sum_{j=1}^dP^S_j(F^{\nu},F^{\nu})\ast^g \Gamma^{v,*}_{\nu,v_j}+\sum_{j=d+1}^{2d}W^S_j(F^{\nu},F^{\nu})\ast^g \Gamma^v_{\nu}
\end{array}
\end{equation} 
We note that the equations in (\ref{prob2aaatrans}) and in (\ref{prob2aaatrans2}) both contain  alternative representations, where the representations in terms of $W^s_j,1\leq j\leq 2d$ are also sufficient for our purposes although these are representations which do not use first order derivatives of the Gaussian or of (the adjoint of) the fundamental solution $\Gamma^v_{\nu}$ respectively. 
Next we consider the functions $W^S_j,~ 1\leq j\leq 2d$, which are related to the representation in (\ref{source1transh}) in more detail. In the following we concentrate on representations in terms of these functions $W^S_j,~ 1\leq j\leq 2d$, where analogous  arguments hold for representations where functions $P^S_j, 1\leq j\leq d$ are involved. We have
\begin{equation}\label{source1transh}
\begin{array}{ll}
Q^S(F,F)=\int_{(w,\sigma)\in {\mathbb R}^3\times S^2}
 {\Big (}\sum_{i=1}^{2d}\int_0^1(1-\theta)\times\\
 \\
\times\partial_{z'_i}G\left( v+\theta \left( -\frac{w}{2}+\frac{\sigma}{2}|w|\right) ,w+2\theta \left( -\frac{w}{2}+\frac{\sigma}{2}|w|\right) \right)d\theta\times\\
\\
\times \left(  -\frac{w_i}{2}+\frac{\sigma_i}{2}|w|\right) {\Big )}|w|dwd\sigma 
=\sum_{j=1}^{2d}W^S_j.
\end{array}
\end{equation}
Note that this representation is with respect $z'=(z'_1,\cdots,z'_{2d})=(v,w)$, where we have introduced $d$ artificial variables $w=(w_1,\cdots ,w_d)$ in the definition of the function $G$ (this function is not to be confused with the Gaussian $G_{\nu}$ of course). We may use approximative representations of the density increment $\delta F=F-F_0$  with an adapted $3d$-dimensional Gaussian $G^{3d}_{\nu}$ (the Gaussian with respect to $(x,z')=(x,v,w)$), and which we denote again by $\delta F^{\nu}$ for simplicity of notation. We have    
\begin{equation}\label{prob2aaaf}
\begin{array}{ll}
\delta F^{\nu}=\sum_{i=1}^dv_iF^{\nu}\ast G^{3d}_{\nu,x_i}
+\sum_{j=1}^{2d}W^S_j(F^{\nu},F^{\nu})\ast G^{3d}_{\nu},
\end{array}
\end{equation} 
where $G^{3d}_{\nu,x_i}$ denotes the spatial derivative of the Gaussian $G^{3d}_{\nu}$ with respect to the variable $x_i$. Note that $\delta F^{\nu}$ depends on the variables $z=(x,v)$ (as the function $G$ is integrated with respect to $w$ in the representation in (\ref{prob2aaaf}). For the first order derivatives with respect to $z_j,~ 1\leq j\leq 2d$ for $z=(x,v)$ we have the representations
\begin{equation}\label{prob2aaafder}
\begin{array}{ll}
\delta F^{\nu}_{,z_j}=F^{\nu}=\sum_{i=1}^dv_iF^{\nu}_{,z_j}\ast G_{\nu,i}
+\sum_{j=1}^{2d}W^S_j(F^{\nu},F^{\nu})\ast G_{\nu,z_j}.
\end{array}
\end{equation} 
We have alternative representations in terms of the fundamental solution $\Gamma^{v,3d}_{\nu}$ of
\begin{equation}
\Gamma^{v,3d}_{\nu,t}-\nu\sum_{j=1}^{3d}\Gamma^{v,3d}_{\nu,\tilde{z}_j,\tilde{z}_j}+\sum_{j=1}^dv_j\Gamma^{v,3d}_{\nu,\tilde{z}_j}=0
\end{equation}  
where $\Gamma^{v,3d}_{\nu}$ depends on $\tilde{z}=(x,v,w)$ including the dummy variables $w=(w_1,\cdots w_d)$. The respective representations are
\begin{equation}\label{prob2aaaf2}
\begin{array}{ll}
\delta F^{\nu}=\sum_{j=1}^{2d}W^S_j(F^{\nu},F^{\nu})\ast^g \Gamma^{v,3d}_{\nu},
\end{array}
\end{equation} 
and
\begin{equation}\label{prob2aaafder}
\begin{array}{ll}
\delta F^{\nu}_{,z_j}=F^{\nu}=\sum_{j=1}^{2d}W^S_j(F^{\nu},F^{\nu})\ast \Gamma^{v,3d}_{\nu,z_j}.
\end{array}
\end{equation} 
In order to obtain an upper bound for all these representations it is essential to estimate second moments of the Gaussian or second moments of Gaussian upper bounds of $\Gamma^v_{\nu}$.Furthemore it is essential to estimate the Gaussian in case of a dimension $d$, and we do this for the sake of notational simplicity.  
Note that the first order spatial derivatives of the Gaussian are antisymmetric in the sense that
\begin{equation}
\left( \frac{-2y_i}{4\nu   t}\right) G_{\nu}(t,y)=\left( \frac{2y^-_i}{4\nu t}\right) G_{\nu}(t,y^-),
\end{equation}
where $y^-=(y^-_1,\cdots,y^-_d)$ with $y^-_i=-y_i$ and $y^-_j=y_j$ for $j\in \lbrace 1,\cdots d\rbrace\setminus \lbrace i\rbrace$.
For any spatially global Lipschitz continuous function $H$ we have
\begin{equation}\label{estest}
\begin{array}{ll}
{\Big |}H\ast G_{\nu,i}=\int H(t-s,x-y)\left( \frac{-2y_i}{4\nu  s}\right) G_{\nu}(s,y)dyds {\Big |}\\
\\
={\Big |} \int \int_{y_i\geq 0} \left( H(t-s,x-y)-H(t-s,x-y^-)\right)\left( \frac{-2y_i}{4\nu  s}\right) G_{\nu}(s,y)dyds{\Big |}\\
\\
\leq {\Big |} L\int \int_{y_i\geq 0} \left( \frac{4y^2_i}{4\nu s}\right) G_{\nu}(s,y)dyds {\Big |}=4LM_2,
\end{array}
\end{equation} 
where $M_2$ is a finite second moment constant of the Gaussian which is small in the following sense. Simplifying a bit more we consider the essential case $d=3$. Note first that for any small time $s>0$ the main mass of the spatial integral is on a ball of radius $\sqrt{ \nu}$. Now, for small $\sqrt{\nu}$ and small time $s$ using partial integration and polar coordinates we observe that (denoting the d-dimensional ball by $B^d_{\sqrt{\nu}}$ and the one dimensional ball by $B_{\sqrt{\nu}}=\left\lbrace r>0|r\leq \sqrt{\nu}\right\rbrace $)  
\begin{equation}\label{series}
\begin{array}{ll}
\frac{1}{4\pi^2}{\Big |} L\int \int_{\left\lbrace y|  y_i\geq 0\right\rbrace \cap B^3_{\sqrt{\nu}}} \left( \frac{4y^2_i}{4\nu   \sigma}\right) G_{\nu}(\sigma,y)dyd\sigma {\Big |}\\
\\
\leq {\Big |} L\int \int_{r\in B_{\sqrt{\nu}}} \left( \frac{4r^2}{4\nu   \sigma}\right) G_{\nu}(\sigma,r)r^2drd\sigma {\Big |}
\\
\\
\leq {\Big |} L\int   \left( \frac{r^5}{5\nu  \sigma}\right) \frac{1}{\sqrt{4\pi \rho\nu \sigma}^D}\exp\left(-\frac{|r|^2}{4\nu \sigma} \right){\Big |}^{\sqrt{\nu}}_{0}d\sigma \\
\\
+ L\int \int_{r\in  B_{\sqrt{\nu}}} \left( \frac{r^5}{5\nu  \rho \sigma}\right) \left( \frac{-2r}{4\nu \sigma}\right)\frac{1}{\sqrt{4\pi \nu \sigma}^D}\exp\left(-\frac{r^2}{4\nu \sigma} \right)drd\sigma {\Big |}\\
\\
\leq {\Big |} L\int   \frac{1}{5\sigma} \frac{1}{\sqrt{4\pi \sigma}^D}\exp\left(-\frac{1}{4\sigma} \right)d\sigma \\
\\
- L\int \int_{r\in  B_{\sqrt{\nu}}} \left( \frac{r^6}{10\nu^2   \sigma^2}\right) \frac{1}{\sqrt{4\pi \nu \sigma}^D}\exp\left(-\frac{r^2}{4\nu \sigma} \right)drd\sigma {\Big |}.
\end{array}
\end{equation}

Iterating this partial integration procedure we get a geometric series upper bound.
More precisely, note that for the last term on the right side in (\ref{series}) we have for $d=3$
\begin{equation}\label{series2}
\begin{array}{ll}
 {\Big |} L\int \int_{r\in  B_{\sqrt{\nu}}} \left( \frac{r^6}{10\nu^2  \sigma^2}\right) \frac{1}{\sqrt{4\pi \nu \sigma}^d}\exp\left(-\frac{r^2}{4\nu \sigma} \right)drd\sigma {\Big |}\\
\\
\leq {\Big |} L\int   \left( \frac{r^7}{70\nu^2  \sigma^2}\right) \frac{1}{\sqrt{4\pi \nu \sigma}^d}\exp\left(-\frac{|r|^2}{4\nu \sigma} \right){\Big |}^{\sqrt{\nu}}_{0}d\sigma {\Big |}\\
\\
+ L\int \int_{r\in  B_{\sqrt{\nu}}} \left( \frac{r^7}{70\nu^2   \sigma^2}\right) \left( \frac{-2r}{4\nu  \sigma}\right)\frac{1}{\sqrt{4\pi \nu \sigma}^d}\exp\left(-\frac{r^2}{4\nu \sigma} \right)drd\sigma {\Big |}\\
\\
\leq {\Big |}L\int   \frac{1}{70\sigma^2} \frac{1}{\sqrt{4\pi \sigma}^D}\exp\left(-\frac{1}{4\sigma} \right)d\sigma \\
\\
- L\int \int_{r\in  B_{\sqrt{\nu}}} \left( \frac{r^8}{140\nu^3   \sigma^3}\right) \frac{1}{\sqrt{4\pi \nu \sigma}^D}\exp\left(-\frac{r^2}{4\nu \sigma} \right)drd\sigma {\Big |}.
\end{array}
\end{equation}
Inductively we have for any $m\geq 1$
\begin{equation}\label{series3}
\begin{array}{ll}
{\Big |} L\int \int_{\left\lbrace y|  y_i\geq 0\right\rbrace \cap B_{\sqrt{\nu}}} \left( \frac{4y^2_i}{4\nu   \sigma}\right) G_{\nu}(\sigma,y)dyd\sigma {\Big |}\\
\\
\leq {\Big |}\sum_{k=1}^m(-1)^{k+1}\rho L\int   \frac{R_k}{\sigma^k} \frac{1}{\sqrt{4\pi \sigma}^D}\exp\left(-\frac{1}{4\sigma} \right)d\sigma {\Big |}\\
\\
+ {\Big |} L\int \int_{r\in  B_{\sqrt{\nu}}} \left( \frac{R_mr^{2m+2}}{2\nu^m   \sigma^m}\right) \frac{1}{\sqrt{4\pi \nu \sigma}^D}\exp\left(-\frac{r^2}{4\nu \sigma} \right)drd\sigma {\Big |},
\end{array}
\end{equation}
where the numbers $R_k$ are recursively defined by
\begin{equation}
R_1=\frac{1}{5},~R_{k+1}=\frac{1}{2}\frac{1}{6+2k-1}R_k,~\mbox{for all $k\geq 1$}.
\end{equation}
%
Clearly,
\begin{equation}
0<R_k<\frac{1}{4^k(k+1)!}.
\end{equation}
We note that
\begin{equation}
\begin{array}{ll}
\frac{d}{d\sigma}\left( \frac{1}{\sigma^{k+d/2}}\right) \exp\left( -\frac{1}{4\sigma}\right)\\
\\
=\left(-\frac{k}{\sigma^{k+1+d/2}}+\frac{1}{4\sigma^{k+2+d/2}} \right) \exp\left( -\frac{1}{4\sigma}\right)\stackrel{!}{=}0\\
\\
\mbox{for}~\sigma=\frac{1}{4k}
\end{array}
\end{equation}
such that the Stirling relation
\begin{equation}
k!\sim \sqrt{2\pi k}\left( \frac{k}{e}\right)^k
\end{equation}
implies that on a small time interval $[0,\Delta]$ we have a converging alternating series on the right sight of (\ref{series3}), i.e. a series which is surely in  $O\left(\frac{1}{\sqrt{k}} \right)$. 
Hence for small $\sqrt{\nu}$ we surely have the upper bound
\begin{equation}
\begin{array}{ll}
{\Big |} L\int \int_{\left\lbrace y| y_i\geq 0\right\rbrace \cap B_{\sqrt{\nu}}} \left( \frac{4y^2_i}{4\nu   \sigma}\right) G_{\nu}(\sigma,y)dyd\sigma {\Big |}\\
\\
\lesssim {\Big |} L\int   \frac{1}{4\sigma} \frac{1}{\sqrt{4\pi \sigma}^d}\exp\left(-\frac{1}{4\sigma} \right)d\sigma {\Big |}.
\end{array}
\end{equation}
Hence, for a small time interval $[t_0,t_0+\Delta_0]$ this upper bound becomes as small as
\begin{equation}\label{parafreeestf}
 {\Big |} L\int_{t_0}^{t_0+\Delta_0}   \frac{1}{4} \frac{1}{\sqrt{4\pi}^D}\frac{1}{(\sigma-t_0)^{2.5}}\exp\left(-\frac{1}{4(\sigma-t_0)} \right)d\sigma {\Big |}\in o(\Delta_0^2),
\end{equation}
where $o$ is the small Landau $o$. Indeed, the left side of (\ref{parafreeestf}) is of polynomial decay, i.e. in $o(\Delta_0^m)$for integer $m\geq 2$, but the statement in (\ref{parafreeestf}) is essentially needed for a global existence scheme based on a time delay transformation which introduces a potential damping term in order to obtain a global regular upper bound which is linear with respect to the time horizon.

Next we construct these global regular upper bounds.  Using the semigroup property it is sufficient to assume that a regular upper bound has been constructed up to a certain time $t_0\geq 0$, and that a linear upper bound can be preserved on a local time interval $[t_0,t_0+\Delta_0]$, where $\Delta_0\in (0,1)$ can be chosen independently of $t_0\geq 0$, but it may depend on the time horizon $\Delta$ 
In order to establish the preservation of a global regular upper bound over a time interval $[t_0,t_0+\Delta_0]$ we use the local analysis above and a comparison function with auto-control.

More precisely, for $t_0\geq 0$ and a time interval $\left[t_{0},t_{0}+\Delta_0 \right]$ for some $\Delta_0 \in (0,1)$ we consider a comparison function $U^{t_0}$ which is constructed by a global time transformation with local time dilatation on a short time interval, i.e., we consider the transformation
\begin{equation}\label{utr1}
\begin{array}{ll}
(1+t )U^{t_{ 0}}(s,y)=F(t,y),\\
\\
\mbox{where}~s=s(t,t_0)=\frac{t-t_{0}}{\sqrt{1-(t-t_{0})^2}},~t-t_{ 0}\in [0,\Delta_0].
\end{array}
\end{equation}

Using the abbreviation $\Delta t=t-t_{0}$, we have
\begin{equation}\label{utr2}
 F_{,t}=
U^{t_{0}}
+(1+t )U^{t_{0}}_{,s}\frac{ds}{dt},~\frac{ds}{dt}
=\frac{1}{\sqrt{1-\Delta t^2}^3}.
\end{equation}
We choose a time interval length $\Delta_0>0$ small enough such that local time contraction holds, and, hence, a local time solution exists on the interval $[t_{ 0},t_{0}+\Delta_0]$. The equation for $U^{t_{0}}$ becomes
\begin{equation}\label{dileq1}
\begin{array}{ll}
U^{t_{ 0}}_{,s}+\sqrt{1-\Delta t^2}^3v\nabla_x U^{t_0}\\
\\
+\frac{\sqrt{1-\Delta t^2}^3}{1+t}Q(t,(1+t )U^{t_{ 0}}(s,.),(1+t )U^{t_{ 0}}(s,.))
=-\frac{\sqrt{1-\Delta t^2}^3}{1+t}U^{t_{0}},\\
\\
~U^{t_{0}}(0,.)=F(t_{ 0},.).
\end{array}
\end{equation}
\begin{rem}
We note that $U^{t_0}$ is a viscosity limit of a sunsequence of the family  $\left( U^{\nu.,t_0}\right)_{\nu >0}$, where $(1+t )U^{\nu,t_{ 0}}(s,y):=F^{\nu}(t,y)$ such that $U^{\nu,t_{ 0}}$ satisfies the equation
\begin{equation}\label{dileq2}
\begin{array}{ll}
U^{\nu,t_{ 0}}_{,s}+\sqrt{1-\Delta t^2}^3 \nu \Delta U^{\nu,t_{0}}+\sqrt{1-\Delta t^2}^3v\nabla_x U^{\nu,t_0}+\\
\\
-\frac{\sqrt{1-\Delta t^2}^3}{1+t}Q(t,(1+t )U^{\nu,t_{ 0}}(s,.),(1+t )U^{\nu,t_{ 0}}(s,.))
=-\frac{\sqrt{1-\Delta t^2}^3}{1+t}U^{\nu,t_{0}},\\
\\
~U^{\nu,t_{0}}(0,.)=F(t_{ 0},.).
\end{array}
\end{equation}
\end{rem}

The equations in (\ref{dileq1}) and (\ref{dileq2}) are defined on a dilated time interval $\left[0,\Delta_{d} \right] $, where $\Delta_d=\frac{\Delta_0}{\sqrt{1-\Delta_0^2}}$.
Now assume that an arbitrary large but fixed time horizon $T>0$ is given. 
Assume that for $m\geq 2$ and for some $0\leq t_0<T$ we have a global upper bound
\begin{equation}
{\Big |}F(t_0,.){\Big |}^{t_0,\Delta_0,s-1}_{sup,1}\leq C(1+t_0)<\infty.
\end{equation}
Here and in the following for time evaluations of functions $F(t_0,.)$ and $U^{t_0}$ we understand that the norm ${\big |}.{\big |}^{t_0,\Delta}_{sup,1}$ is evaluated for that fixed time.
Then 
\begin{equation}
{\Big |}U^{t_{ 0}}(0,.){\Big |}^{t_0,\Delta_0,s-1}_{sup,1}=\frac{1}{1+t_{ 0}} {\Big |}F(t_{0},.){\Big |}^{t_0,\Delta_0,s-1}_{sup,1}\leq C.
\end{equation}
Then the local analysis above shows that for $\Delta_0\leq \frac{1}{(1+T)^2}$ small enough we have
\begin{equation}
\begin{array}{ll}
{\Big |}\frac{\sqrt{1-\Delta t^2}^3}{1+t}Q(t,(1+t )U^{t_{ 0}}(s,.),(1+t )U^{t_{ 0}}(s,.)){\Big|}^{t_0,\Delta_0,s-1}_{sup,1}
\leq \frac{C}{2}\Delta_0.
\end{array}
\end{equation}
We may assume that $C\geq 1$ and studying two cases where $${\Big |}U^{t_{ 0}}{\Big |}^{t_0,\Delta_0,s-1}_{sup,1}\leq \frac{C}{2},{\Big |} U^{t_{ 0}}{\Big |}^{t_0,\Delta_0,s-1}_{sup,1}\in \left[ \frac{C}{2},C\right] $$ we conclude that for $\Delta_0\leq \frac{1}{(1+T)^2}$ small enough we have
\begin{equation}
{\Big |}U^{t_{ 0}}(\Delta_d,.){\Big |}^{t_0,\Delta_0,s-1}_{sup,1}\leq C.
\end{equation} 
This implies that 
\begin{equation}
\begin{array}{ll}
{\Big |}F(t_{0}+\Delta_0,.){\Big |}^{t_0,\Delta_0,s-1}_{sup,1}
\leq (1+t_{ 0}+\Delta_0){\Big |}U^{\nu,t_{ 0}}(\Delta_d,.){\Big |}^{t_0,\Delta_0,s-1}_{sup,1}\\
\\
 \leq C(1+t_{ 0}+\Delta_0).
\end{array}
\end{equation} 
Iterating this procedure for all integers $k\geq 0$ with $k\Delta_0 \leq T$ we have inductively
\begin{equation}
\begin{array}{ll}
{\Big |}F(k\Delta_0,.){\Big |}^{(k-1)\Delta_0,\Delta_0,s-1}_{sup,1}\leq {\Big |}U^{(k-1)\Delta_0}(0,.){\Big |}^{(k-1)\Delta_0,\Delta_0,s-1}_{sup,1}\leq C(1+k\Delta_0).
\end{array}
\end{equation}
Applying the local time once again we can interpolate and obtain for all $0\leq t\leq T$
\begin{equation}
\begin{array}{ll}
{\Big |}F(t,.){\Big |}^{0,T,s-1}_{sup,1}\leq (1+t){\Big |}U^{0}_i(0,.){\Big |}^{0,\Delta_0,s-1}_{sup,1}
\leq C(1+t)
\end{array}
\end{equation}
by inductive reasoning.

The same reasoning as in this item iii) holds for approximating densities $F^{\nu},~ \nu>0$ and for the sequence $(F^{\nu_k})_k$ of the viscosity limit construction with $\nu_k\downarrow 0$, of course. For the latter series $F^{\nu_k}\geq 0$ holds for all $\nu_k$ and also in the viscosity limit  $\nu_k\downarrow 0$. Here we have also an alternative reasoning compared to that in  (\cite{PL}). 
Assuming that $F^{\nu_k}(t_0,.)\geq 0$ and ${\big |}F^{\nu_k}(t_0,.){\big |}_{L^1}> 0$ we 
may compare $\nu_k(\Delta F^{\nu_k}(t_0,.))\ast G^D_{\nu}$ and $F^{\nu_k}\ast G_{\nu_k,t}$ for some dimension $D\geq 3$ and are a hypothetical minimum $(t_m,x_m)$ for small $\nu_k$ and use our local analysis of the nonlinear terms above in order to show that $F^{\nu_k}$ stays nonnegative over time.

%
%
%
%
\end{itemize}

\appendix

\section{Levy expansion of $\Gamma^{\nu}_v$}
A priori estimates and smoothness of densities for linear second order parabolic equation with smooth and linearly bounded coefficients are well-known from stochastic analysis. However, some footnotes and alternative arguments concerning the spacial circumstances of the linear part of the Boltzmann equation and concerning the adjoint may be useful.  The main theorem and its proof can be found in \cite{KS}. We recall
\begin{thm}
	\label{stroock}
	Consider a $d$-dimensional difffusion process of the form
	\begin{equation}\label{stochm}
		\mathrm{d}X_t \ = \ \sum_{i=1}^d\sigma_{0i}(X_t)\mathrm{d}t + \sum_{j=1}^{d}\sigma_{ij}(X_t)\mathrm{d}W^j_t
	\end{equation}
with $X(0)=x\in {\mathbb R}^d$ with values in ${\mathbb R}^d$ and on a time interval $[0,T]$. Here $W_t$ denotes a standard $d$-dimensional Brownian motion, and we assume that $\sigma_{0i},\sigma_{ij}\in C^{\infty}_{lb}$, where $C^{\infty}_{lb}$ denotes the space of linearly bounded smooth functions with bounded derivatives.  Then the law of the process $X$ is absolutely continuous with respect to the Lebesgue measure, and the density $p$ exists and is smooth, i.e. 
\begin{equation}
\begin{array}{ll}
p:(0,T]\times {\mathbb R}^d\times{\mathbb R}^d\rightarrow {\mathbb R}\in C^{\infty}\left( (0,T]\times {\mathbb R}^d\times{\mathbb R}^d\right). 
\end{array}
\end{equation}
Moreover, for each nonnegative natural number $j$, and multiindices $\alpha,\beta$ there are increasing functions of time
\begin{equation}\label{constAB}
A_{j,\alpha,\beta}, B_{j,\alpha,\beta}:[0,T]\rightarrow {\mathbb R},
\end{equation}
and functions
\begin{equation}\label{constmn}
n_{j,\alpha,\beta}, 
m_{j,\alpha,\beta}:
{\mathbb N}\times {\mathbb N}^d\times {\mathbb N}^d\rightarrow {\mathbb N},
\end{equation}
such that 
\begin{equation}\label{pxest}
\begin{array}{ll}
{\Bigg |}\frac{\partial^j}{\partial t^j} \frac{\partial^{|\alpha|}}{\partial x^{\alpha}} \frac{\partial^{|\beta|}}{\partial y^{\beta}}p(t,x,y){\Bigg |}\\
\\
\leq \frac{A_{j,\alpha,\beta}(t)(1+x)^{m_{j,\alpha,\beta}}}{t^{n_{j,\alpha,\beta}}}\exp\left(-B_{j,\alpha,\beta}(t)\frac{(x-y)^2}{t}\right)
\end{array}
\end{equation}
Moreover, all functions (\ref{constAB}) and  (\ref{constmn}) depend on the level of iteration of Lie-bracket iteration at which the H\"{o}rmander condition becomes true.
\end{thm}

The alternative preceding argument involving $\Gamma^{\nu}_v$  requires the existence of this density in $C^1$. We point out that this follows form classical analysis via the Levy expansion. 
Formally the Levy expansion of the fundamental solution $\Gamma^{\nu}_v$ of
\begin{equation}\label{gausextapp}
\partial_t \Gamma^v_{\nu}+\nu\Delta \Gamma^v_{\nu}+v\nabla_x\Gamma^v_{\nu}=0
\end{equation}
is of the form
\begin{equation}
\begin{array}{ll}
\Gamma^{\nu}_v(t,x,v;s,y)=G_{\nu}(t,x,v;s,y)\\
\\
+\sum_{k=1}^{\infty}\int_0^t\int_{{\mathbb R}^{2d}}G_{\nu}(t,x,v;\sigma,z)L^{G_{\nu}}_{FO_k}(\sigma,z,s,y)dzds,~z,y\in {\mathbb R}^{2d},
\end{array}
\end{equation}
where $L^{G_{\nu}}_{FO_k}$ are convolutions of the first order operator part in (\ref{gausextapp}) defined recursively by
\begin{equation}
\begin{array}{ll}
L^{G_{\nu}}_{FO_1}:=v\nabla_xG_{\nu},~
L^{G_{\nu}}_{FO_{k+1}}:=L^{G_{\nu}}_{FO_1}\ast^g L^{G_{\nu}}_{FO_k},~ k\geq 1,
\end{array}
\end{equation}
where
\begin{equation}
\begin{array}{ll}
\left(  L^{G_{\nu}}_{FO_1}\ast^g L^{G_{\nu}}_{FO_k}\right) (t,x;s,y)=\int_0^t\int_{{\mathbb R}^{2d}}L^{G_{\nu}}_{FO_1}(t,x,v;\sigma,z)L^{G_{\nu}}_{FO_k}(\sigma,z,s,y)dzds
\end{array}
\end{equation}
It follows from standard arguments in the case of bounded coefficients slightly adapted to this situation with linear drift coefficients that this prescription defines the fundamental solution for all  $(x,v)\in {\mathbb R}^{2d}$ and $t>s> 0$. Moreover, we observe that we get a linear factor in front of the Gaussian a priori upper bounds of $\Gamma^v_{\nu}$ and its first derivatives with respect to spatial variables due to the linear coefficients $v_i$.
Finally a remark concerning the adjoint.

 Greens's identity implies that the fundamental solution $\Gamma$ of a linear parabolic equation
of second order, where
\begin{equation}
\Gamma_{,t}-\sum_{ij}^da_{ij}\Gamma_{,i,j}-\sum_j^d b_j\Gamma_{,j}=0
\end{equation} 
with regular linearly bounded variable coefficients, and the usual uniform ellipticity condition satisfies the condition
\begin{equation}\label{adj}
D^{\alpha}_z\Gamma(t,x;s,y)=D^{\alpha}_z\Gamma^*(s,y;t,x),
\end{equation}
where $\Gamma*$ is thesolution of the adjoint equation, $|\alpha|\geq 0$, and the derivatives $D^{\alpha}_z$ refer to the first spatial variables ($x=(x_1,\cdots,x_d)$ on the right side and $y=(y_1,\cdots,y_d)$ on the left side.  Furthermore the upper bound for $|\alpha|$ is only limited by the joint regularity of $\Gamma$ and $\Gamma^*$. This relation in (\ref{adj}) can be used in situations where a standard convolution rule does not hold.

\end{document}